\documentclass[12pt]{amsart}
\usepackage{amsmath,amssymb,amsthm,amsrefs}
\usepackage{enumitem}
\usepackage{color,hyperref}
\usepackage{color}
\hypersetup{colorlinks,breaklinks,
	linkcolor=blue,urlcolor=blue,
	anchorcolor=blue,citecolor=blue}
\usepackage{cleveref}
\theoremstyle{plain}
\newtheorem{theorem}{Theorem}[section]
\newtheorem{proposition}[theorem]{Proposition}
\newtheorem{corollary}[theorem]{Corollary}

\theoremstyle{definition}
\newtheorem{definition}[theorem]{Definition}

\newtheorem{remark}[theorem]{Remark}

\makeatletter
\@addtoreset{equation}{section}
\makeatother

\makeatletter
\newcommand{\Spvek}[2][r]{%
  \gdef\@VORNE{1}
  \left(\hskip-\arraycolsep%
    \begin{array}{#1}\vekSp@lten{#2}\end{array}%
  \hskip-\arraycolsep\right)}

\def\vekSp@lten#1{\xvekSp@lten#1;vekL@stLine;}
\def\vekL@stLine{vekL@stLine}
\def\xvekSp@lten#1;{\def\temp{#1}%
  \ifx\temp\vekL@stLine
  \else
    \ifnum\@VORNE=1\gdef\@VORNE{0}
    \else\@arraycr\fi%
    #1%
    \expandafter\xvekSp@lten
  \fi}
\makeatother

\begin{document}
\title[Spectral bounds for periodic Jacobi matrices]
{Spectral bounds for periodic Jacobi matrices\\}
\author{Burak Hat\.{i}no\u{g}lu}
\address{Department of Mathematics, Michigan State University, East Lansing MI 48829, U.S.A.}
\email{hatinogl@msu.edu}




\begin{abstract}
We consider periodic Jacobi operators and obtain upper and lower estimates on the sizes of the spectral bands. Our proofs are based on estimates on the logarithmic capacities and connections between the Chebyshev polynomials and logarithmic capacity of compact subsets of the real line. 
\end{abstract}
\maketitle

\section{Introduction}\label{Sec1}

The periodic Jacobi operator $J$ with period $p$, acting on the Hilbert space $l^2(\mathbb{Z})$ is the self-adjoint operator associated with the infinite Jacobi matrix
\begin{equation*}
 \begin{pmatrix}
\ddots & \ddots & 0 & 0 & 0 \\
\ddots & b_1 & a_1 & 0 & 0 \\
0 & a_1 & b_2 & a_2 & 0 \\
0 & 0  &  a_2 & b_3 & \ddots  \\
0 & 0 & 0 & \ddots & \ddots
\end{pmatrix},
\end{equation*}
given by
\begin{equation*}
    (J\psi)_n = a_{n-1}\psi_{n-1} + b_n\psi_n + a_n\psi_{n+1}, 
\end{equation*}
where $a_n \in \mathbb{R}$, $b_n > 0$, $a_{n+p} = a_n$ and $b_{n+p} = b_n$ for any $n \in \mathbb{Z}$. The periodic discrete Schr\"{o}dinger operators on the real line give a subclass, where $a_n = 1$ for any $n \in \mathbb{Z}$. It is well known that the spectrum $\sigma$ of the periodic Jacobi operator $J$ with period $p \geq 2$ is purely absolutely continuous and consists of $p$ possibly touching (but not overlapping) bands on the real line
$$
\sigma_n = [\lambda_n^{min},\lambda_n^{max}]
$$
for $1 \leq n \leq p$, seperated by the spectral gaps
$$
\gamma_n = (\lambda_n^{max},\lambda_{n+1}^{min})
$$
for $1 \leq n \leq p-1$ \cite{vM76,T89,T00}. 

In this paper, we discuss bounds for the spectral bands and gaps of periodic Jacobi operators and obtain new results using logarithmic potential theory. 

The endpoints of the spectrum are bounded by Gershgorin circle theorem as
$$
\min_{n} (b_n - a_n - a_{n-1}) \leq \lambda_1^{min} \qquad \text{and} \qquad \lambda_p^{max} \leq \max_{n} (b_n + a_n + a_{n-1}).
$$
The distance between maximal and minimal points of the spectrum is also estimated below
\begin{equation}\label{est1}
    s := \lambda_p^{max} - \lambda_1^{min} \geq 4A,
\end{equation}
where $A := (a_1a_2\cdots a_p)^{1/p}$ \cite{KK03}. The same estimate is an upper bound for the Lebesgue measure of the spectrum
\begin{equation}\label{est2}
  |\sigma| = \sum_{n=1}^p |\sigma_n| \leq 4A,  
\end{equation}
where $|\cdot|$ denotes the Lebesgue measure \cite{KK03,L92}. Another upper estimate for the Lebesgue measure of the spectrum is given in terms of diagonal entries of the Jacobi matrix and the distance between maximal and minimal points of the spectrum 
\begin{equation}\label{est3}
     \sum_{n=1}^p |\sigma_n| \leq s - m,
\end{equation}
where $s = \lambda_p^{max} - \lambda_1^{min}$ and $m = \max_{n} b_n - \min_{n} b_n$ \cite{DS83,L92}.
The upper estimate \eqref{est2} was improved to
\begin{equation}\label{est4}
 \sum_{n=1}^p |\sigma_n| \leq 4\min_{n} a_n   
\end{equation}
in \cite{K15}. The Lebesgue measure of the spectrum is also estimated below as
\begin{equation}\label{est5}
     \sum_{n=1}^p |\sigma_n| \geq \frac{4A^p}{M^{p-1}},
\end{equation}
where $M = \max\big\{\max_{n}(b_n + a_n + a_{n-1}) - \min_{n}b_n, \max_{n}b_n - \min_{n}(b_n - a_n - a_{n-1})\big\}$ \cite{DS83,L92}, but $M$ can be replaced by $s$ \cite{KK03}, so we have the lower estimate
\begin{equation}\label{est6}
     \sum_{n=1}^p |\sigma_n| \geq \frac{4A^p}{s^{p-1}}.
\end{equation}\label{est7}
Combining \eqref{est1} with \eqref{est4} gives a lower estimate for the total size of the gaps
\begin{equation}\label{est8}
    \sum_{n=1}^{p-1}|\gamma_n| \geq 4(A - \min_{n}a_n).
\end{equation}
A lower estimate in terms of diagonal entries of the Jacobi matrix is a direct consequence of \eqref{est3} and the band-gap structure of the spectrum:
\begin{equation}\label{est9}
     \sum_{n=1}^{p-1}|\gamma_n| \geq \max_{n} b_n - \min_{n} b_n.
\end{equation}
Another lower estimate of the total size of the gaps is
\begin{equation}
    \sum_{n=1}^{p-1}|\gamma_n| \geq \max\Big\{\max\{4A,2\max_{n}a_n\} - 4\min_{n}a_n~,~ \max_{n}b_n - \min_{n}b_n\Big\},
\end{equation}
which was obtained in \cite{K19}.

Some similar bounds were obtained for quasi-periodic \cite{PR09} and matrix-valued \cite{K15} Jacobi operators.

We obtain the following estimates for the spectral bands using logarithmic potential theory and Chebyshev polynomials. 

\begin{theorem}\label{LogSumLower}
    Let $J$ be a periodic Jacobi operator with period $p$ defined by off-diagonal and diagonal sequences $\{a_n\}_{n \in \mathbb{Z}}$ and $\{b_n\}_{n \in \mathbb{Z}}$, respectively. Let $s:= \lambda_p^{max} - \lambda_1^{min}$, $A := (a_1a_2\cdots a_p)^{1/p}$ and $\sigma_n$ denote the nth spectral band of $J$. If $s \leq d$, then
    \begin{equation}
    \frac{1}{\log\big(d/A\big)}\leq\displaystyle\sum_{n=1}^{p}\frac{1}{\log\big(4d/|\sigma_{n}|\big)}.
    \end{equation}
\end{theorem}

\begin{corollary}\label{LogSingleLower}
    Let $J$ be a periodic Jacobi operator with period $p$ defined by off-diagonal and diagonal sequences $\{a_n\}_{n \in \mathbb{Z}}$ and $\{b_n\}_{n \in \mathbb{Z}}$, respectively. Let $s:= \lambda_p^{max} - \lambda_1^{min}$, $A := (a_1a_2\cdots a_p)^{1/p}$ and $\sigma_n$ denote the nth spectral band of $J$. Then
    \begin{equation}
    \frac{4A^p}{s^{p-1}} \leq \max_{1 \leq n \leq p} |\sigma_n|.
    \end{equation}
\end{corollary}

\begin{theorem}\label{LogSumUpper}
    Let $J$ be a periodic Jacobi operator with period $p$ defined by off-diagonal and diagonal sequences $\{a_n\}_{n \in \mathbb{Z}}$ and $\{b_n\}_{n \in \mathbb{Z}}$, respectively. Let $A := (a_1a_2\cdots a_p)^{1/p}$, and $\sigma_n$ and $\gamma_n$ denote the nth spectral band and nth spectral gap of $J$, respectively. If $$\displaystyle\min_{1 \leq n \leq p-1}|\gamma_n| \geq d,$$ then
    \begin{equation}
    \frac{1}{\log^+\big(d/A\big)}\geq\displaystyle\sum_{n=1}^{p}\frac{1}{\log^+\big(4d/|\sigma_{n}|\big)}.
    \end{equation}
\end{theorem}

\begin{corollary}\label{LogSingleUpper}
    Let $J$ be a periodic Jacobi operator with period $p$ defined by off-diagonal and diagonal sequences $\{a_n\}_{n \in \mathbb{Z}}$ and $\{b_n\}_{n \in \mathbb{Z}}$, respectively. Let $A := (a_1a_2\cdots a_p)^{1/p}$, and $\sigma_n$ and $\gamma_n$ denote the nth spectral band and gap of $J$, respectively. If the condition 
    $$\displaystyle 4\Big(\min_{1 \leq n \leq p-1}|\gamma_n|\Big) \geq \max\Big\{\max_{1 \leq n \leq p}|\sigma_n|~,~ 4A\Big\}$$
    is satisfied, then
    \begin{equation}
     \min_{1 \leq n \leq p} |\sigma_n| \leq \frac{4A^p}{\Big(\displaystyle\min_{1 \leq n \leq p-1}|\gamma_n|\Big)^{p-1}}.
    \end{equation}
\end{corollary}

The paper is organized as follows.
 
Section \ref{Sec2} includes preliminaries required for our proofs, namely some basics of logarithmic potential theory, definition and properties of Chebyshev polynomials, and the discriminant of a periodic Jacobi operator and its fundamental properties.

Section \ref{Sec3} includes proofs of Theorems \ref{LogSumLower} and \ref{LogSumUpper}, and Corollaries \ref{LogSingleLower} and \ref{LogSingleUpper}.

\section{Preliminaries}\label{Sec2}

In our proofs we use estimates on logarithmic capacity, so let's recall some basics of logarithmic potential theory, which can be found e.g. in \cite{R95}.

\begin{definition}
Let $\mu$ be a finite Borel measure, supported on a compact subset of the complex plane. Then \textit{logarithmic potential} of $\mu$ is the function $U^{\mu}:\mathbb{C}\rightarrow(-\infty,\infty]$ defined by 
\begin{equation*}
U^{\mu}(z):=\int\log\frac{1}{|z-\omega|}d\mu(\omega)
\end{equation*}
\end{definition}

\begin{definition}
Let $\mu$ be a finite Borel measure, supported on a compact subset of the complex plane. Its \textit{logarithmic energy} $I(\mu)\in(-\infty,\infty]$ is defined by
\begin{equation*}
I(\mu):=\int U^{\mu}(z)d\mu(z)=\int\int\log\frac{1}{|z-\omega|}d\mu(\omega)d\mu(z)
\end{equation*}
\end{definition}

\begin{definition}
Let $K$ be a compact subset of $\mathbb{C}$ and $M(K)$ be the set of Borel probability measures compactly supported inside $K$. The measure $\mu_{K}\in M(K)$ is called the \textit{equilibrium measure} for $K$ if $$I(\mu_{K})=\displaystyle\inf_{\mu\in M(K)}I(\mu).$$
\end{definition}

Now we are ready to define the logarithmic capacity.

\begin{definition}
The \textit{logarithmic capacity} of a subset $E$ of the complex plane is given by
\begin{equation*}
Cap(E):=\sup_{\mu\in M(E)}\exp(-I(\mu)).
\end{equation*}
In particular if $K$ is compact with equilibrium measure $\mu_{K}$, then $Cap(K)=\exp(-I(\mu_{K}))$.
\end{definition}

Logarithmic capacity is monotone and multiplicative with respect to modulus, but it is not weakly subadditive. However, in terms of its logarithm some estimations can be made on logarithmic capacity of union.

\begin{theorem} \label{union cap}\emph{\textbf{(}\cite{R95}, \textit{Theorem 5.1.4}\textbf{)}}
Let $E:=\cup_{n=1}^{m}E_{n}$ be the union of Borel subsets $E_{n}$ of $\mathbb{C}$ and $d > 0$, where $m\in\mathbb{N}\cup\{\infty\}$.
\begin{enumerate}
    \item If $diam(E)\leq d$, then
\begin{equation*}
\displaystyle\frac{1}{\log(\frac{d}{Cap(E)})}\leq\sum_{n=1}^{m}\frac{1}{\log(\frac{d}{Cap(E_{n})})}
\end{equation*}
	 \item If $dist(E_{j},E_{k})\geq d$ whenever $j\neq k$, then
\begin{equation*}
\frac{1}{\log^{+}(\frac{d}{Cap(E)})}\geq\sum_{n=1}^{m}\frac{1}{\log^{+}(\frac{d}{Cap(E_{n})})}
\end{equation*}
\end{enumerate}
Note that in the first part the logarithms are always nonnegative since $Cap(E_n) \leq Cap(E) \leq diam(E) \leq d$.
\end{theorem}

Theorem \ref{union cap} will be the main tool in our proofs. We will also consider logarithmic capacity of the spectrum and the spectral bands. Logarithmic capacity of an interval is one fourth of its length. In order to understand logarithmic capacity of the spectrum, we need to recall relations between logarithmic capacity and Chebyshev polynomials.

\begin{definition}
Let $K$ be a compact subset of $\mathbb{C}$ and $T_{n,K}$ be the monic polynomial of degree $n$ such that $||{T}_{n,K}||_{K}\leq||P||_{K}$ for any monic polynomial $P$ of degree $n$, where $||\cdot||_{K}$ denotes the uniform norm over $K$. Then $T_{n,K}$ is called the \emph{nth Chebyshev polynomial on $K$} and $||{T}_{n,K}||_{K}$ is called the \textit{nth Chebyshev number of K}, denoted by $t_{n}(K)$.
\end{definition}

The alternation theorem will allow us to make the connection with the spectrum and its Chebyshev polynomial.

\begin{definition}\label{altset}
    Let $P$ be a real polynomial of degree $n$. Then $P$ has an \textit{alternating set} in $K \subset \mathbb{R}$ if there exists $\{x_k\}_{k=0}^n$ in $K$ satisfying $x_0 < x_1 < \cdots < x_n$ such that
$$
P(x_k) = (-1)^{n-k}||P||_{K}.
$$
\end{definition}

\begin{theorem}\label{alternation}\normalfont \textbf{(Alternation Theorem)} \emph{\textbf{(}\cite{CSZ17}, \textit{Theorem 1.1}\textbf{)}}
    Let $K$ be a compact subset of the real line. The nth Chebyshev polynomial on $K$ has an alternating set in $K$. Conversely, any monic polynomial with an alternating set in $K$ is the Chebyshev polynomial on $K$.
\end{theorem}

The sequence of Chebyshev numbers of a compact set is not necessarily convergent. However subadditivity of logarithms of Chebyshev numbers imply existence of $\lim_{n \rightarrow \infty} t_{n}^{1/n}(K)$. This limit is called the Chebyshev number of $K$ and by a classical result of Szeg\"{o} it is nothing but the logarithmic capacity of $K$. The following result of Peherstorfer and Totik will help us to understand the logarithmic capacity of the spectrum. 

\begin{theorem}\label{TotikPeherstorfer} \emph{\textbf{(}\cite{T11}, \textit{Theorem 1} \textrm{ and } \cite{P11}, \textit{Proposition 1.1}}\textbf{)}
Let $K=\displaystyle\cup_{j=1}^{l}[a_{j},b_{j}]$. Also let $T_{n,K}$ and $t_{n}(K)$ denote the nth Chebyshev polynomial and nth Chebyshev number of $K$, respectively. For a natural number
$n\geq 1$ the following are pairwise equivalent.
\begin{itemize}
\item[a)] $\displaystyle\frac{t_n(K)}{{Cap(K)}^n}=2$.\\

\item[b)] $T_{n,K}$ has $n+l$ extreme points on $K$.\\

\item[c)] $K=\{z$ $|$ $T_{n;K}(z)\in[-t_{n}(K),t_{n}(K)]\}$.\\

\item[d)] If $\mu_{K}$ denotes the equilibrium measure of $K$, then for each $j=1,2,\dots,l$,
$$
\mu_{K}([a_{j},b_{j}]) = \frac{q_j}{n},
$$
where $q_{j}+1$ is the number of extreme points on $[a_{j},b_{j}]$.\\

\item[e)] With $\pi(x)=\prod_{j=1}^{l}(x-a_{j})(x-b_{j})$ the equation
\begin{equation*}
P_{n}^{2}(x)-\pi(x)Q_{n-l}^{2}(x)= c
\end{equation*}
is solvable for the polynomials $P_{n}$ and $Q_{n-l}$ of degree $n$ and $n-l$, respectively, where $c$ is a positive constant.
\end{itemize}
\end{theorem}

\begin{remark}
The ratio in item (a) is called Widom Factor \cite{GH15}, denoted by $W_{n}(K)$, and its estimates and asymptotics are studied recently \cite{AGH16,CSZ17,CSZ22}.  
\end{remark}

Next, we need to recall the discriminant of a periodic Jacobi operator. If $p$ is the period of the Jacobi operator $J$, then its discriminant is defined by
\begin{equation}
    \Delta(\lambda) := \textrm{tr} \Big(\prod_{n=p}^{1} A_n(\lambda)\Big),
\end{equation}
where 
$$
\displaystyle A_n(\lambda) := \begin{pmatrix}
\frac{\lambda - b_n}{a_n} & -\frac{a_{n-1}}{a_n}\\
1 & 0
\end{pmatrix}.
$$
The discriminant is the trace of the monodromy matrix, which is defined as the matrix that shifts by the period of the Jacobi operator $J$ along the solutions of the eigenvalue problem $J\psi = \lambda\psi$, i.e.
$$
\begin{pmatrix}
\psi_{n+1}\\
\psi_{n}
\end{pmatrix} = A_n(\lambda) \begin{pmatrix}
\psi_{n}\\
\psi_{n-1}
\end{pmatrix},
$$
if $\psi = \{\psi_n\}_{n \in \mathbb{Z}}$ is a solution of $J\psi = \lambda\psi$.

Floquet-Bloch theory allows us to get the following properties of the discriminant $\Delta(\lambda)$ and its connection with the spectrum of the corresponding periodic Jacobi operator \cite{T89,T00}.

\begin{proposition}\label{PropertiesofDiscriminant}
    Let $\Delta$ be the discriminant of a Jacobi operator $J$ with period $p$. Then we have the following properties. 
    \begin{enumerate}
        \item The discriminant $\Delta(\lambda)$ is a polynomial in the spectral parameter $z$ of degree $p$.
        \item The leading coefficient of the discriminant $\Delta(\lambda)$ is $1/(a_1a_2\cdots a_p)$.
        \item The discriminant $\Delta(\lambda)$ splits over the real line with $p$ distinct roots.
        \item If $\lambda_0$ is a local maxima or minima of $\Delta(\lambda)$, then $|\Delta(\lambda_0)| \geq 2$.
        \item The spectrum $\sigma(J)$ is the inverse image of $[-2,2]$ under $\Delta(\lambda)$ on the real line, consists of $p$ possibly touching but not overlapping bands and is purely absolutely continuous.
    \end{enumerate}
\end{proposition}

\section{Proofs}\label{Sec3}

\begin{proof}[\normalfont \textbf{Proof of Theorem~\ref{LogSumLower}}] 
Since each spectral band $\sigma_n$ is an interval, we have
$$
Cap(\sigma_n) = \frac{|\sigma_n|}{4}.
$$
In order to get logarithmic capacity of the spectrum $\sigma$ of $J$, using items 1 and 2 of Proposition \ref{PropertiesofDiscriminant} recall that the discriminant $\Delta(E)$ of $J$ is a polynomial of degree $p$ with the leading coefficient $1/A^p := 1/(a_1a_2\cdots a_p)$, so $A^p\Delta(E)$ is a monic polynomial. Using items 3,4 and 5 of Proposition \ref{PropertiesofDiscriminant} we observe that the discriminant $\Delta$ has an alternating set in the spectrum $\sigma$ according to Definition \ref{altset}. Therefore by Alternation Theorem, Theorem \ref{alternation}, the discriminant is the pth Chebyshev polynomial of the spectrum, i.e.
$$
T_{p,\sigma}(E) = \Delta(E)
$$
and
$$
t_p(\sigma) = 2A^p.
$$
Note that by item 5 of Proposition \ref{PropertiesofDiscriminant}, the spectrum satisfies item c of Theorem \ref{TotikPeherstorfer}, so using item a of the same theorem we get $t_p(\sigma)/Cap(\sigma)^p = 2$, and hence the logarithmic capacity of the spectrum as
$$
Cap(\sigma) = A.
$$
Now recalling that $s = diam(\sigma)$ and $s \leq d$, and using the first part of Theorem \ref{union cap} we get the desired result
 \begin{equation*}
    \frac{1}{\log\big(d/A\big)}\leq\displaystyle\sum_{n=1}^{p}\frac{1}{\log\big(4d/|\sigma_{n}|\big)}.
    \end{equation*}
\end{proof}

\begin{proof}[\normalfont \textbf{Proof of Corollary~\ref{LogSingleLower}}] 
Using Theorem \ref{LogSumLower} and noting that $|\sigma_n| < s \leq d$ for any $1 \leq n \leq p$ we get 
\begin{equation}\label{maxineq}
    \frac{1}{\log\big(d/A\big)}\leq\displaystyle\sum_{n=1}^{p}\frac{1}{\log\big(4d/|\sigma_{n}|\big)} \leq \frac{p}{\log\big(4d/\max_{1 \leq n \leq p} |\sigma_n|\big)}. 
    \end{equation}
We also know that $d/A > 1$, so logarithms on both ends of \eqref{maxineq}  are positive. Therefore replacing $d$ with $s$ we get
$$
\frac{4s}{\displaystyle\max_{1 \leq n \leq p} |\sigma_n|} \leq \Big(\frac{s}{A}\Big)^p,
$$
and hence the desired result
$$
\frac{4A^p}{s^{p-1}} \leq \displaystyle\max_{1 \leq n \leq p} |\sigma_n|.
$$
\end{proof}

\begin{proof}[\normalfont \textbf{Proof of Theorem~\ref{LogSumUpper}}] 
Since each spectral band $\sigma_n$ is an interval
$$
Cap(\sigma_n) = \frac{|\sigma_n|}{4}.
$$
In order to get logarithmic capacity of the spectrum $\sigma$ of $J$, using items 1 and 2 of Proposition \ref{PropertiesofDiscriminant} recall that the discriminant $\Delta(E)$ of $J$ is a polynomial of degree $p$ with the leading coefficient $1/A^p := 1/(a_1a_2\cdots a_p)$, so $A^p\Delta(E)$ is a monic polynomial. Using items 3,4 and 5 of Proposition \ref{PropertiesofDiscriminant} we observe that the discriminant $\Delta$ has an alternating set in the spectrum $\sigma$ according to Definition \ref{altset}. Therefore by Alternation Theorem, Theorem \ref{alternation}, the discriminant is the pth Chebyshev polynomial of the spectrum, i.e.
$$
T_{p,\sigma}(E) = \Delta(E)
$$
and
$$
t_p(\sigma) = 2A^p.
$$
Note that by item 5 of Proposition \ref{PropertiesofDiscriminant}, the spectrum satisfies item c of Theorem \ref{TotikPeherstorfer}, so using item a of the same theorem we get $t_p(\sigma)/Cap(\sigma)^p = 2$, and hence the logarithmic capacity of the spectrum as
$$
Cap(\sigma) = A.
$$
Since the spectrum has band-gap structure and we assume $\min_{1 \leq n \leq p-1}|\gamma_n| \geq d$, we get $dist(\sigma_{j},\sigma_{k})\geq d$ whenever $j \neq k$. Therefore using the second part of Theorem \ref{union cap} we get the desired result
    \begin{equation*}
    \frac{1}{\log^+\big(d/A\big)}\geq\displaystyle\sum_{n=1}^{p}\frac{1}{\log^+\big(4d/|\sigma_{n}|\big)}.
    \end{equation*}
\end{proof}

\begin{proof}[\normalfont \textbf{Proof of Corollary~\ref{LogSingleUpper}}] 
The condition we assumed allow us to replace positive logarithms in Theorem \ref{LogSumUpper} with logarithms, and $d$ with $\min_{1 \leq n \leq p-1}|\gamma_n|$, so we have
    \begin{align*}
    \frac{1}{\log\big(\min_{1 \leq n \leq p-1}|\gamma_n|/A\big)} &\geq \displaystyle\sum_{j=1}^{p}\frac{1}{\log\big(4\min_{1 \leq n \leq p-1}|\gamma_n|/|\sigma_{j}|\big)}\\ &\geq \frac{p}{\log\big(4\min_{1 \leq n \leq p-1}|\gamma_n|/\min_{1 \leq n \leq p} |\sigma_n|\big)}.
    \end{align*}
Therefore we get
$$
4\frac{\displaystyle\min_{1 \leq n \leq p-1}|\gamma_n|}{\displaystyle\min_{1 \leq n \leq p} |\sigma_n|} \geq \Big(\frac{\displaystyle\min_{1 \leq n \leq p-1}|\gamma_n|}{A}\Big)^p,
$$
and hence the desired result
$$
\displaystyle\min_{1 \leq n \leq p} |\sigma_n| \leq \frac{4A^p}{\displaystyle\Big(\min_{1 \leq n \leq p-1}|\gamma_n|\Big)^{p-1}}.
$$
\end{proof}

\section{Acknowledgments}
Part of this work was conducted at Georgia Institute of Technology, where the author was a postdoc of Svetlana Jitomirskaya. The author thanks funding from NSF DMS-2052899, DMS-2155211, and Simons 681675.

\bibliographystyle{abbrv}
\bibliography{references}

\end{document}